\DeclareSymbolFont{AMSb}{U}{msb}{m}{n}
\documentclass[10pt,a4paper,leqno,noamsfonts]{amsart}
   \makeatletter
   \renewcommand\@biblabel[1]{#1.}
   \makeatother
\usepackage[english]{babel}
\usepackage[dvipsnames]{xcolor}
\usepackage{graphicx,pifont,soul,physics} 
\usepackage[utopia]{mathdesign}
\usepackage[a4paper,top=3cm,bottom=3cm,
           left=3cm,right=3cm,marginparwidth=60pt]{geometry}
\usepackage[utf8]{inputenc}
\usepackage{braket,caption,comment,mathtools,stmaryrd}
\usepackage[usestackEOL]{stackengine}
\usepackage{multirow,booktabs,microtype,relsize}
\usepackage[foot]{amsaddr}
\usepackage[colorlinks,bookmarks]{hyperref} %
      \hypersetup{colorlinks,%
            citecolor=britishracinggreen,%
            filecolor=black,%
            linkcolor=cobalt,%
            urlcolor=cornellred}
      \numberwithin{equation}{section}
\usepackage[capitalise]{cleveref}
\usepackage[colorinlistoftodos]{todonotes}
\definecolor{antiquewhite}{rgb}{0.98, 0.92, 0.84}
\definecolor{buff}{rgb}{0.94, 0.86, 0.51}
\definecolor{palecopper}{rgb}{0.85, 0.54, 0.4}
\definecolor{fluorescentyellow}{rgb}{0.8, 1.0, 0.0}
\definecolor{bole}{rgb}{0.47, 0.27, 0.23}

\usepackage{amsmath,float}
\usetikzlibrary{decorations.markings}

\definecolor{cornellred}{rgb}{0.7, 0.11, 0.11}
\definecolor{britishracinggreen}{rgb}{0.0, 0.26, 0.15}
\definecolor{cobalt}{rgb}{0.0, 0.28, 0.67}
\DeclareSymbolFont{usualmathcal}{OMS}{cmsy}{m}{n}
\DeclareSymbolFontAlphabet{\mathcal}{usualmathcal}

\newcommand{\TT}{\mathbf{T}}

\newcommand{\BA}{{\mathbb{A}}}

\newcommand{\BC}{{\mathbb{C}}}

\newcommand{\BG}{{\mathbb{G}}}

\newcommand{\BL}{{\mathbb{L}}}

\newcommand{\BN}{{\mathbb{N}}}

\newcommand{\BP}{{\mathbb{P}}}

\newcommand{\BZ}{{\mathbb{Z}}}

\newcommand{\CI}{{\mathcal I}}

\newcommand{\CT}{{\mathcal T}}

\newcommand{\simto}{\,\widetilde{\to}\,}

\newcommand{\fix}{\mathsf{fix}}

\DeclareMathOperator{\Hilb}{Hilb}
\DeclareMathOperator{\Sets}{Sets}
\DeclareMathOperator{\Sch}{Sch}
\DeclareMathOperator{\Quot}{Quot}

\DeclareMathOperator{\coker}{coker}

\DeclareMathOperator{\Sym}{Sym}

\DeclareMathOperator{\Coh}{Coh}

\DeclareMathOperator{\Exp}{Exp}
\DeclareMathOperator{\Var}{Var}

\newcommand{\into}{\hookrightarrow}
\newcommand{\onto}{\twoheadrightarrow}

\DeclareFontFamily{OT1}{rsfs}{}
\DeclareFontShape{OT1}{rsfs}{n}{it}{<-> rsfs10}{}
\DeclareMathAlphabet{\curly}{OT1}{rsfs}{n}{it}
\renewcommand\hom{\mathscr{H}\kern-0.3em\mathit{om}}

\newcommand\Ext{\operatorname{Ext}}
\newcommand\Hom{\operatorname{Hom}}

\DeclareMathOperator{\lHom}{\mathscr{H}\kern-0.3em\mathit{om}}
\DeclareMathOperator{\RRlHom}{\mathbf{R}\kern-0.025em\mathscr{H}\kern-0.3em\mathit{om}}

\DeclareMathOperator{\lExt}{{\mathscr{E}\kern-0.2em\mathit{xt}}}

\DeclareMathOperator{\pr}{pr}

\newcommand\Spec{\operatorname{Spec}}
\newcommand\Supp{\operatorname{Supp}}

\newcommand\id{\operatorname{id}}

\newcommand{\OO}{\mathscr O}



\usepackage[all]{xy}
\usepackage{tikz}
\usepackage{tikz-cd}
\usepackage{adjustbox}
\usepackage{rotating}
\usepackage{comment}

\usetikzlibrary{matrix,shapes,intersections,arrows,decorations.pathmorphing}
\tikzset{commutative diagrams/arrow style=math font}
\tikzset{commutative diagrams/.cd,
mysymbol/.style={start anchor=center,end anchor=center,draw=none}}
\newcommand\MySymb[2][\square]{%
  \arrow[mysymbol]{#2}[description]{#1}}
\tikzset{
shift up/.style={
to path={([yshift=#1]\tikztostart.east) -- ([yshift=#1]\tikztotarget.west) \tikztonodes}
}
}

\theoremstyle{definition}

\newtheorem*{lemma*}{Lemma}
\newtheorem*{theorem*}{Theorem}
\newtheorem*{example*}{Example}
\newtheorem*{fact*}{Fact}
\newtheorem*{notation*}{Notation}
\newtheorem*{definition*}{Definition}
\newtheorem*{prop*}{Proposition}
\newtheorem*{remark*}{Remark}
\newtheorem*{corollary*}{Corollary}

\newtheorem*{conventions*}{Conventions}

\newtheorem{definition}{Definition}[section]

\newtheorem{example}[definition]{Example}

\newtheorem{notation}[definition]{Notation}
\newtheorem{remark}[definition]{Remark}

\newtheoremstyle{thm} 
        {3mm}
        {3mm}
        {\slshape}
        {0mm}
        {\bfseries}
        {.}
        {1mm}
        {}
\theoremstyle{thm}
\newtheorem{theorem}[definition]{Theorem}
\newtheorem{corollary}[definition]{Corollary}
\newtheorem{lemma}[definition]{Lemma}
\newtheorem{prop}[definition]{Proposition}
\newtheorem{thm}{Theorem}

\newtheoremstyle{ex} 
        {3mm}
        {3mm}
        {}
        {0mm}
        {\scshape}
        {.}
        {1mm}
        {}
\theoremstyle{ex}

\newtheoremstyle{sol} 
        {3mm}
        {3mm}
        {}
        {0mm}
        {\scshape}
        {.}
        {1mm}
        {}
\theoremstyle{sol}

\newtheorem*{Acknowledgments*}{Acknowledgments}

\newcommand{\boldit}[1]{\boldsymbol{#1}} 
\newcommand{\bfk}{\mathbf{k}}
\newcommand{\nn}{\boldit{n}}

\newcommand{\Filt}{\mathrm{Filt}}
\DeclareMathAlphabet\BCal{OMS}{cmsy}{b}{n}

\title[On the motive of the nested Quot scheme of points on a curve]{On the motive of the nested Quot scheme of points on a curve}

\author{Sergej Monavari}
\address{Mathematical Institute, Utrecht University, P.O.~Box 80010 3508 TA Utrecht, The Netherlands}
\email{s.monavari@uu.nl}
\author{Andrea T. Ricolfi}
\address{SISSA Trieste, Via Bonomea 265, 34136 Trieste, Italy}
\email{aricolfi@sissa.it}

\date{}

\begin{document}
\maketitle

\begin{abstract}
Let $C$ be a smooth curve over an algebraically closed field $\bfk$, and let $E$ be a locally free sheaf of rank $r$. We compute, for every $d>0$, the generating function of the motives $[\Quot_C(E,\nn)] \in K_0(\Var_{\bfk})$, varying $\nn = (0\leq n_1\leq\cdots\leq n_d)$, where $\Quot_C(E,\nn)$ is the \emph{nested Quot scheme of points}, parametrising $0$-dimensional subsequent quotients $E \onto T_d \onto \cdots \onto T_1$ of fixed length $n_i = \chi(T_i)$. The resulting series, obtained by exploiting the Białynicki-Birula decomposition, factors into a product of shifted motivic zeta functions of $C$. In particular, it is a rational function. 
\end{abstract}


\section{Introduction}

Let $K_0(\Var_{\bfk})$ be the Grothendieck ring of varieties over an algebraically closed field $\bfk$. If $Y$ is a $\bfk$-variety, its \emph{motivic zeta function} 
\[
\zeta_Y(q) = 1 + \sum_{n>0}\,\bigl[\Sym^nY\bigr] q^n\,\in\,K_0(\Var_{\bfk})\llbracket q \rrbracket
\]
is a generating series introduced by Kapranov in \cite{Kapranov_rational_zeta}, where he proved that for smooth curves it is a rational function in $q$.

In this paper we compute the motive of the \emph{nested Quot scheme of points} $\Quot_C(E,\nn)$ on a smooth curve $C$, entirely in terms of $\zeta_C(q)$. Here, $E$ is a locally free sheaf on $C$, and $\nn = (0 \leq n_1\leq\cdots\leq n_d)$ is a non-decreasing tuple of integers, for some fixed $d>0$. The scheme $\Quot_C(E,\nn)$ generalises the classical Quot scheme  of Grothendieck (recovered when $d=1$): it parametrises flags of quotients $E \onto T_d \onto \cdots \onto T_1$ where $T_i$ is a $0$-dimensional sheaf of length $n_i$.

Our main result, proved in \Cref{mainthm_BODY} in the main body, is the following. 

\begin{thm}\label{thm:main_thm_INTRO}
Let $C$ be a smooth curve over $\bfk$, let $E$ be a locally free sheaf of rank $r$ on $C$. Then
\[
\sum_{0\leq n_1\leq\cdots\leq n_d}\,\bigl[\Quot_C(E,\nn)\bigr]q_1^{n_1}\cdots q_d^{n_d} = \prod_{\alpha=1}^r\prod_{i=1}^d\zeta_C\left(\BL^{\alpha-1}q_iq_{i+1}\cdots q_d\right)\,\in\,K_0(\Var_{\bfk})\llbracket q_1,\ldots,q_d \rrbracket,
\]
where $\BL=[\BA^1_{\bfk}]$ is the Lefschetz motive.
In particular, this generating function is rational in $q_1,\ldots,q_d$.
\end{thm}

The statement taken with $d=1$, thus regarding the motive $[\Quot_C(E,n)]$ of the usual Quot scheme of points, was proved in \cite{BFP19}. Our result is a natural generalisation, which was inspired by Mochizuki's paper on ``Filt schemes'' \cite{Mochizuki_FILT}.

\smallbreak
Our formula fits nicely in the philosophical path according to which
\begin{center}
``rank $r$ theories factorise in $r$ rank $1$ theories''.
\end{center}
There are to date a number of examples of this phenomenon in Donaldson--Thomas theory, exhibiting a generating series of rank $r$ invariants as a product of $r$ (suitably shifted) generating series of rank $1$ invariants: see for instance \cite{BR18,Virtual_Quot} for enumerative DT invariants, \cite{FMR_K-DT} for K-theoretic DT invariants, \cite{cazzaniga2020higher,Cazzaniga:2020aa} for motivic DT invariants and \cite{Magnificent_colors,dZNPZ_playing_index_M_theory} for the parallel pictures in string theory.

\smallbreak

\smallbreak
The paper is organised as follows.
In Section \ref{sec: nested quot} we introduce the \emph{nested Quot scheme} and prove its connectedness. In Section \ref{sec: tangent space} we describe its tangent space and prove that, for a smooth quasiprojective curve, the nested Quot scheme is smooth. Under the assumption that the locally free sheaf is split, in \Cref{sec:BB} we describe a torus framing action and its associated  Białynicki-Birula decomposition.  In Section \ref{sec: motives} we prove that the motive of the  nested Quot scheme is independent of the locally free sheaf, and exploit the Białynicki-Birula decomposition to prove \Cref{thm:main_thm_INTRO}. Our result readily implies closed formulae for the generating series of Hodge--Deligne polynomials, $\chi_y$-genera, Poincar\'e polynomials, Euler characteristics, since these are all motivic measures; we provide some explicit formulae in \Cref{HD_poly}.\\

After our paper was written, we were informed that our formula for the motive of the nested Quot scheme on a \emph{projective} curve can be alternatively obtained, after some manipulations, from general results on the stack of iterated Hecke correspondences \cite[Corollary 4.10]{Heinloth_altri} (see also \cite[Section~3]{hoskins2019formula} for a related computation of the Voevodsky motive with rational coefficients). Our paper provides a direct and self-contained argument for this formula, exploiting the geometry of the nested Quot scheme.

\subsection*{Acknowledgments}
We thank Alina Marian for spotting a mistake in the first draft of this paper. We also wish to thank Joachim Jelisiejew for stimulating conversations around the Białynicki-Birula decomposition. We are grateful to Takuro Mochizuki for helpful discussions about nested Quot schemes. We thank Jochen Heinloth, Victoria Hoskins and Simon Pepin Lehalleur for carefully explaining the relation to their previous work. 
The first author is grateful to Martijn Kool for useful discussions. We finally thank the anonymous referee for several helpful comments, which improved the exposition of the paper.

S.M.~is supported by NWO grant TOP2.17.004. A.R.~is funded by Dipartimenti di Eccellenza.

\begin{conventions*}
All \emph{schemes} are of finite type over an algebraically closed field $\bfk$. A \emph{variety} is a reduced separated $\bfk$-scheme. If $Y$ is a scheme and $Y_1,\ldots,Y_s$ are locally closed subschemes of $Y$, we say that they form a (locally closed) \emph{stratification}, denoted `$Y=Y_1\amalg\cdots\amalg Y_s$' with a slight abuse of notation, if the natural morphism of schemes $Y_1\amalg\cdots\amalg Y_s \to Y$ is bijective. This is crucial in our calculations since this condition implies the identity $[Y] = [Y_1]+\cdots+[Y_s]$ in $K_0(\Var_{\bfk})$.
\end{conventions*}

\section{Nested Quot schemes of points}\label{sec: nested quot}

\subsection{The moduli space}
Let $X$ be a quasiprojective $\bfk$-variety and $E$ a coherent sheaf on $X$. Fix an integer $d>0$ and a non-decreasing $d$-tuple $\nn=(n_1\leq \dots \leq n_d)$ of non-negative integers $n_i\in \BZ_{\geq 0}$. We define the \emph{nested Quot functor} associated to $(X,E,\nn)$ to be the functor $\mathsf{Quot}_X(E,\nn) \colon \Sch_{\bfk}^{\mathrm{op}} \to \Sets$ sending a $\bfk$-scheme $B$ to the set of isomorphism classes of subsequent quotients
\begin{align*}
    E_B\twoheadrightarrow \CT_d \twoheadrightarrow \dots \twoheadrightarrow \CT_1,
\end{align*}
where $E_B$ is the pullback of $E$ along $X\times_\bfk B \to X$ and $\CT_i \in \Coh(X\times_\bfk B)$ is a $B$-flat family of $0$-dimensional sheaves of length $n_i$ over $X$ for all $i=1,\ldots,d$. Two `nested quotients' 
\begin{align*}
    E_B\twoheadrightarrow \CT_d \twoheadrightarrow \dots \twoheadrightarrow \CT_1, \qquad E_B\twoheadrightarrow \CT_d' \twoheadrightarrow \dots \twoheadrightarrow \CT_1'
\end{align*}
are considered isomorphic when $\ker (E_B \onto \CT_i) = \ker (E_B \onto \CT_i')$ for all $i=1,\ldots,d$.

The representability of the functor $\mathsf{Quot}_X(E,\nn)$ can be proved adapting the proof of \cite[Theorem 4.5.1]{Sernesi} or by an explicit induction on $d$ as in \cite[Section 2.A.1]{modulisheaves}. We define $\Quot_X(E, \nn)$ to be the moduli scheme representing the above functor. Its closed points are then in bijection with the set of isomorphism classes of nested quotients
\begin{align*}
    E\twoheadrightarrow T_d \twoheadrightarrow \dots \twoheadrightarrow T_1,
\end{align*}
where each $T_i \in \Coh(X)$ is a $0$-dimensional quotient of $E$ of length $n_i$. 
The nested Quot scheme comes with a closed immersion
\begin{equation}\label{eqn:embedding_in_product}
\Quot_X(E,\nn) \into \prod_{i=1}^d \Quot_X(E,n_i)
\end{equation}
cut out by the nesting condition $\ker (E\onto T_d) \into \ker(E\onto T_{d-1}) \into \cdots\into \ker (E\onto T_1)$. In particular, it is projective as soon as $X$ is projective. If $C$ is a smooth proper curve over $\BC$ and $E \in \Coh(C)$ is a locally free sheaf, the cohomology of $\Quot_C(E,\nn)$ was studied by Mochizuki \cite{Mochizuki_FILT}.

\begin{example}
The classical Quot scheme $\Quot_X(E,n)$ of length $n$ quotients of $E$ is obtained by setting $\nn = (n)$, i.e.~taking $d=1$ and $n_d=n$. If we set $\nn = (1\leq 2\leq \cdots\leq d)$, we obtain Mochizuki's \emph{complete Filt scheme} $\Filt(E,d)$, which for $d=1$ reduce to $\Filt(E,1)=\BP(E)$ \cite{Mochizuki_FILT}.  When $E=\OO_X$, we use the notation $\Hilb^{\nn}(X)$ to denote $\Quot_X(\OO_X,\nn)$. This space is the \emph{nested Hilbert scheme of points}, studied extensively by Cheah \cite{MR1616606,MR2716513,MR1612677}.
\end{example}

\subsection{Support map and nested punctual Quot scheme}
Fix a variety $X$, a coherent sheaf $E$ and a $d$-tuple of non-negative integers $\nn = (n_1\leq \cdots\leq n_d)$ for some $d>0$. Composing the embedding \eqref{eqn:embedding_in_product} with the usual Quot-to-Chow morphisms yields the \emph{support map}
\begin{equation}\label{supportmap}
\mathsf h_{E,\nn}\colon \Quot_X(E,\nn) \into \prod_{i=1}^d \Quot_X(E,n_i) \to \prod_{i=1}^d \Sym^{n_i}(X)
\end{equation}
recording the $0$-cycles $([\Supp T_i]\in \Sym^{n_i}(X))_{1\leq i\leq d}$ attached to a $d$-tuple $(E \onto T_i)_{1\leq i\leq d}$. Here, $\Sym^mX = X^m/\mathfrak S_m$ is the $m$-th symmetric power of $X$.

We make the following definition.

\begin{definition}[Nested punctual Quot scheme]
Let $X$ be a variety, $x\in X$ a point, $E\in \Coh(X)$ a coherent sheaf, $\nn = (n_1\leq \cdots\leq n_d)$ a tuple of non-negative integers. The \emph{nested punctual Quot scheme} attached to $(X,E,\nn,x)$ is the closed subscheme
\[
\Quot_X(E,\nn)_x \subset \Quot_X(E,\nn),
\]
defined as the preimage of the cycle $(n_1x,\ldots,n_dx)$ along the support map $\mathsf h_{E,\nn}$.
\end{definition}

The name `punctual' refers, as for the classical Quot schemes, to the fact that all quotients are entirely supported at a single point. We will not need the following result.

\begin{lemma}
Let $X$ be a smooth quasiprojective variety of dimension $m$, and let $E$ be a locally free sheaf of rank $r$ on $X$. For every $d$-tuple $\nn = (n_1\leq \cdots\leq n_d)$, and for every $x \in X$, one has a non-canonical isomorphism
\[
\Quot_X(E,\nn)_x \cong \Quot_{\BA^m}(\OO^{\oplus r},\nn)_0.
\]
\end{lemma}

\begin{proof}
The result follows from the isomorphism $\Quot_X(E,k)_x\simto\Quot_{\BA^m}(\OO^{\oplus r},k)_0$ relating the classical punctual Quot schemes, which is proved in full detail in \cite[Section 2.1]{ricolfi2019motive} exploiting a choice of \'etale coordinates around $x$ (which exist by the smoothness assumption, and which explain the non-canonical nature of the isomorphism). It remains to observe that the induced isomorphism
\[
\begin{tikzcd}
\displaystyle\prod_{i=1}^d \Quot_X(E,n_i)_x \arrow{r}{\sim} & \displaystyle\prod_{i=1}^d \Quot_{\BA^m}(\OO^{\oplus r},n_i)_0 
\end{tikzcd}
\]
maps the subscheme $\Quot_X(E,\nn)_x$ isomorphically onto $\Quot_{\BA^m}(\OO^{\oplus r},\nn)_0$.
\end{proof}

\subsection{Connectedness}
We prove the following connectedness result for the nested Quot scheme. A proof in the case $(r,d,\nn)=(1,1,n)$ of the classical Hilbert scheme was first given by Hartshorne \cite{MR0213368}, and by Fogarty in the surface case \cite{Fogarty_Hilb}. We shall also exploit Cheah's connectedness result for $\Hilb^{\nn}(X)$, see~\cite[Sec. 0.4]{MR1616606}.

\begin{theorem}\label{thm: connected} If $X$ is an irreducible quasiprojective $\bfk$-variety and $E$ is a locally free sheaf on $X$, then $\Quot_X(E,\nn)$ is connected for every $\nn = (n_1\leq \cdots \leq n_d)$. In particular, the classical Quot scheme $\Quot_X(E,n)$ is connected for every $n\geq 0$.
\end{theorem}

\begin{proof}
The proof consists of several steps.

\textsc{Step 1}: We reduce to proving the statement when $E = \OO_X^{\oplus r}$ is trivial. Let $x = [E \onto T_d\onto \cdots\onto T_1] \in \Quot_X(E,\nn)$ be a point, where $E$ is arbitrary. Since $T_d$ is $0$-dimensional we can find an open neighbourhood $U \subset X$ of the set-theoretic support of $T_d$ such that $E|_U = \OO_U^{\oplus r}$ is trivial. The point $x$ then lies in the image of the open immersion $\Quot_U(\OO_U^{\oplus r},\nn) \into \Quot_X(E,\nn)$. By assumption, the space $\Quot_U(\OO_U^{\oplus r},\nn)$ is connected. Now if $x' = [E \onto T_d'\onto \cdots\onto T_1'] \in \Quot_X(E,\nn)$ is another point, we can find another open subset $U' \subset X$ surrounding the support of $T'_d$ and trivialising $E$. Since $X$ is irreducible, $U \cap U' \neq \emptyset$, which implies $\Quot_U(\OO_U^{\oplus r},\nn) \cap \Quot_{U'}(\OO_{U'}^{\oplus r},\nn) \neq \emptyset$, so $x$ and $x'$ are connected in $\Quot_X(E,\nn)$ by any point in this intersection.


\textsc{Step 2}: The scheme $\Quot_X(\OO_X^{\oplus r},\nn)$ has a framing $\TT$-action with non-empty fixed locus, where $\TT = \BG_m^r$ (see \Cref{prop:fixed_locus} for an explicit description of this fixed locus: we shall exploit it in the next step). Let $x \in \Quot_X(\OO_X^{\oplus r},\nn)$ be an arbitrary point. Then the closure of its orbit contains a $\TT$-fixed point --- this will be explained in \Cref{sec:BB}. Therefore it is enough to prove that any two $\TT$-fixed points $x,x' \in \Quot_X(\OO_X^{\oplus r},\nn)^\TT$ are connected in $\Quot_X(\OO_X^{\oplus r},\nn)$.

\textsc{Step 3}: In principle, we should check connectedness for an \emph{arbitrary} pair $(x,x')$ of $\TT$-fixed points
\begin{align*}
x &= [\OO_X^{\oplus r} \onto T_d\onto \cdots\onto T_1]  \,\,\in\,\, \prod_{\alpha=1}^r \Hilb^{\nn_\alpha}(X)\subset \Quot_X(\OO_X^{\oplus r},\nn)^\TT,
\\
x' &= [\OO_X^{\oplus r} \onto T_d'\onto \cdots\onto T_1']  \,\,\in\,\, \prod_{\alpha=1}^r \Hilb^{\nn'_\alpha}(X)\subset \Quot_X(\OO_X^{\oplus r},\nn)^\TT,
\end{align*}
where $\sum_{1\leq\alpha\leq r}\nn_\alpha=\nn=\sum_{1\leq\alpha\leq r}\nn_\alpha'$. But since each nested Hilbert scheme $\Hilb^{\boldit{m}}(X)$ is connected (cf.~\cite[Sec. 0.4]{MR1616606}), we can in fact choose a pair of convenient $x$ and $x'$. We fix them satisfying the condition that $\Supp(T_d),\Supp(T_d')$  consist of $n_d$ distinct points. When viewed in the full space $\Quot_X(\OO_X^{\oplus r},\nn)$, the points $x$ and $x'$ both belong to the open subset 
\[
U \subset \Quot_X(\OO_X^{\oplus r},\nn),
\]
defined by the cartesian diagram
\begin{equation}\label{eqn: diagramma U}
\begin{tikzcd}[row sep=large]
U\MySymb{dr}\arrow[r]    \arrow[d, hook]& \prod_{i=1}^d(\Sym^{n_i}X\setminus \Delta_{\textrm{big}})\arrow[d, hook,"\mathrm{open}"]\\
\Quot_X(\OO_X^{\oplus r},\nn)\arrow{r}{\mathsf h_{\OO_X^{\oplus r},\nn}} & \prod_{i=1}^d\Sym^{n_i}X
\end{tikzcd}
\end{equation}
where $\Delta_{\textrm{big}}\subset \Sym^{n_i}X$ is the big diagonal and the bottom map is the support map \eqref{supportmap}.
In other words, $U\subset \Quot_X(\OO_X^{\oplus r},\nn) $ is the open subscheme consisting of the flags of quotients $ [\OO_X^{\oplus r}\twoheadrightarrow T_d\twoheadrightarrow\dots \twoheadrightarrow T_1]$ where each $T_i$ is supported on $n_i$ distinct points.  This yields an open immersion
\[
U\into \prod_{i=1}^d V_i,
\]
where $V_i\subset \Quot_X(\OO_X^{\oplus r}, n_i-n_{i-1})$ is the open subscheme consisting of points $[\OO_X^{\oplus r}\twoheadrightarrow T'_i]$ where the quotients $T'_i$ are supported on $n_i-n_{i-1}$ distinct points (and we set $n_0=0$). 
The scheme $V_i$ is the image of the \'etale map (cf.~\cite[Proposition A.3]{BR18})
\[
\begin{tikzcd}
 A_i \arrow{r}{\oplus} & \Quot_X(\OO_X^{\oplus r},n_i-n_{i-1})
\end{tikzcd}
\]
defined on the open subscheme
\[
A_i \subset \Quot_X(\OO_X^{\oplus r}, 1)^{n_i-n_{i-1}}
\]
parametrising quotients $(\OO_X^{\oplus r} \onto \OO_{x_k})_{k}$ with $x_k \neq x_l$ for every $k\neq l$. On the other hand, 
\[
\Quot_X(\OO_X^{\oplus r}, 1)^{n_i-n_{i-1}} \cong \BP(\OO_X^{\oplus r})^{n_i-n_{i-1}} \cong (X \times_{\bfk} \BP^{r-1})^{n_i-n_{i-1}}
\]
is irreducible, hence $A_i$ is irreducible, and in particular $V_i$ is irreducible, being the image of an irreducible space along a continuous map. Therefore $U \into \prod_iV_i$ is also irreducible, in particular connected, which completes the proof.
\end{proof}

\section{Tangent space and smoothness in the case of curves}\label{sec: tangent space}

Fix $(X,E,\nn)$ as in the previous section. 
For any point $x\in \Quot_X(E, \nn)$ representing a $d$-tuple of nested quotients
\[
\begin{tikzcd}
E \arrow[two heads]{r} &
T_d \arrow[two heads]{r}{p_{d-1}} &
T_{d-1} \arrow[two heads]{r}{p_{d-2}} &
\cdots \arrow[two heads]{r}{p_2} &
T_2 \arrow[two heads]{r}{p_{1}} &
T_1 
\end{tikzcd}
\]
we set $K_i = \ker (E \onto T_i)$. We have a flag of subsheaves
\[
\begin{tikzcd}
K_d \arrow[hook]{r}{\iota_{d-1}} & 
K_{d-1} \arrow[hook]{r}{\iota_{d-2}} & 
\cdots\arrow[hook]{r}{\iota_2} &
K_2 \arrow[hook]{r}{\iota_1} &
K_1 \arrow[hook]{r} & E
\end{tikzcd}
\]
and, for any $i=1,\ldots,d-1$, maps
\begin{align*}
    \phi_i\colon \Hom_X(K_i, T_i)&\to \Hom_X(K_{i+1}, T_i), & g \mapsto g\circ \iota_{i} \\
    \psi_i\colon \Hom_X(K_{i+1}, T_{i+1})&\to \Hom_X(K_{i+1}, T_{i}), & h \mapsto p_i\circ h 
\end{align*}
which we assemble in a  matrix
\[
\Delta_x = 
\begin{pmatrix}
-\phi_1 & \psi_1 & 0 & 0 & \cdots & 0 \\
0 & -\phi_2 & \psi_2 & 0 & \cdots & 0 \\
\vdots & \vdots & \vdots & \vdots & \vdots & \vdots \\
0 & 0 & 0 & \cdots & -\phi_{d-1} & \psi_{d-1}
\end{pmatrix}
\]
defining a map
\[ 
\begin{tikzcd}
\Delta_x\colon \displaystyle\bigoplus_{i=1}^d\Hom_X(K_i, T_i)
\arrow{r} & 
\displaystyle\bigoplus_{i=1}^{d-1}\Hom_X(K_{i+1}, T_{i}).
\end{tikzcd}
\]
The embedding \eqref{eqn:embedding_in_product} induces a $\bfk$-linear inclusion of tangent spaces
\[
T_x\Quot_X(E,\nn) \into \bigoplus_{i=1}^d \Hom_X(K_i,T_i),
\]
which can be described as follows: a $d$-tuple of maps $(\delta_1,\ldots,\delta_d) \in \bigoplus_{i=1}^d \Hom_X(K_i,T_i)$ belongs to the tangent space of $\Quot_X(E,\nn)$ at $x$ precisely when the diagram
\begin{equation}\label{eqn: comm diagram tangent space}
\begin{tikzcd}[row sep=large,column sep=large]
  K_d\arrow[hook]{r}{\iota_{d-1}} \arrow{d}{\delta_d}
& K_{d-1} \arrow[hook]{r}{\iota_{d-2}} \arrow{d}{\delta_{d-1}}
& \cdots \arrow[hook]{r}{\iota_2}
& K_2\arrow[hook]{r}{\iota_1}\arrow{d}{\delta_2}
& K_1 \arrow{d}{\delta_1} \\
  T_d \arrow[two heads]{r}{p_{d-1}} 
& T_{d-1} \arrow[two heads]{r}{p_{d-2}}
& \cdots \arrow[two heads]{r}{p_2}
& T_2 \arrow[two heads]{r}{p_1}
& T_1
\end{tikzcd}
\end{equation}

commutes. This is formalised in terms of the map $\Delta_x$ in the next proposition.

\begin{prop}
Set $\nn=(n_1\leq \dots \leq n_d)$. The tangent space of $\Quot_X(E, \nn)$ at a point $x=[E\twoheadrightarrow T_d\twoheadrightarrow\dots \twoheadrightarrow T_1]$ is
\begin{align*}
   T_x\Quot_X(E, \nn)=\ker\left(\bigoplus_{i=1}^d\Hom(K_i, T_i)\xrightarrow{\Delta_x} \bigoplus_{i=1}^{d-1}\Hom(K_{i+1}, T_{i})\right).
\end{align*}
In particular, if $E$ is locally free of rank $r$ on a smooth curve $C$, we have that $\Quot_C(E, \nn)$ is smooth of dimension $r \cdot n_d$.
\end{prop}
\begin{proof}
Along the same lines of \cite[Prop. 4.5.3(i)]{Sernesi} it is easy to see that the tangent space is given by the maps making Diagram \eqref{eqn: comm diagram tangent space} commute, which is equivalent to belonging to the kernel of $\Delta_x$.

Let $Q_i$ be the 0-dimensional sheaf fitting in the exact sequences 
\begin{align*}
   & 0\to K_i\to K_{i-1}\to Q_i\to 0\\
&0\to Q_i \to T_i\to T_{i-1}\to 0
\end{align*}
for every $i=1,\dots, d$.
If $X=C$ is a smooth curve, we have that each $K_i$ is a locally free sheaf of rank $r$ (because torsion free is equivalent to locally free on smooth curves); since $Q_i$ is a $0$-dimensional sheaf, we obtain the vanishings
\begin{align}\label{eqn: vanishings}
\Ext^j_C(K_i, T_i)=\Ext^j_C(K_{i+1}, T_i)=\Ext^j_C(K_i, Q_i)=0, \quad j>0.
\end{align}
Therefore each $\psi_i$ is a surjective map, which implies that $\Delta_x$ is surjective and that the dimension of the tangent space is  computed as 
\begin{align*}
    \dim_\bfk T_x\Quot_C(E, \nn)&=\dim_\bfk\left( \bigoplus_{i=1}^d\Hom_C(K_i, T_i)\right)- \dim_\bfk\left(\bigoplus_{i=1}^{d-1}\Hom_C(K_{i+1}, T_{i})\right) \\
   &= \sum_{i=1}^d rn_i-\sum_{i=1}^{d-1} rn_i\\
    &= rn_d.
\end{align*}
Since the tangent space dimension is constant and $\Quot_C(E,\nn)$ is connected by \Cref{thm: connected}, it is enough to find a smooth open subset $U \subset \Quot_C(E,\nn)$ of dimension $r n_d$. We shall exploit the fact that the classical Quot scheme $\Quot_C(E,m)$ is smooth of dimension $rm$, which follows from standard deformation theory and the vanishing $\Ext_C^1(K,T) = \mathrm{H}^1(C,K^\vee\otimes T) = 0$ for an arbitrary point $[K \into E \onto T]\in\Quot_C(E,m)$.

Let $U\subset\Quot_C(E,\nn)$ be the open subscheme as in Diagram \eqref{eqn: diagramma U} (which of course exists for arbitrary $E$), and write $U\cong \prod_{i=1}^dV_i$ as in the proof of   \Cref{thm: connected}. We know that each $V_i\subset \Quot_C(E, n_i-n_{i-1})$ is smooth of dimension $r\cdot (n_i-n_{i-1})$, therefore $U$ is smooth of dimension $r n_d$ as required.
\end{proof}
\begin{remark}
The smoothness of $\Quot_C(E, \nn)$ was already proved by Mochizuki \cite[Prop.~2.1]{Mochizuki_FILT}, via a tangent-obstruction theory argument. See also  \cite{Lissite-quot} for the classification of smoothness of $ \Quot_X(E, \nn)$ when $X$ has arbitrary dimension.
\end{remark}

\section{Białynicki-Birula decomposition}\label{sec:BB}
Let $E$ be a locally free sheaf of rank $r$ on a variety $X$. Assume that $E=\bigoplus_{\alpha=1}^r L_\alpha$ splits into a sum of line bundles on $X$. Then $\Quot_X(E,\nn)$ admits the action of the algebraic torus $\TT=\BG_m^r$ as in \cite{Bifet}. Indeed, $\TT$ acts diagonally on the product $\prod_{i=1}^d \Quot_X(E,n_i)$ and the closed subscheme $\Quot_X(E,\nn)$ is $\TT$-invariant. Its fixed locus is determined by a straightforward generalisation of the main result of \cite{Bifet}.

\begin{prop}\label{prop:fixed_locus}
If $E=\bigoplus_{\alpha=1}^r L_\alpha$, there is a scheme-theoretic identity
\begin{align*}
    \Quot_X\left(E, \nn\right)^{\TT} = \coprod_{\substack{\nn_1 + \cdots + \nn_r = \nn }}\prod_{\alpha=1}^r   \Quot_X(L_\alpha, \nn_\alpha).
\end{align*}
\end{prop}
\begin{proof}
We construct a bijection on $\bfk$-valued points, which is straightforward to verify in families.

Fix tuples $\nn_\alpha = (n_{\alpha,1}\leq \dots \leq n_{\alpha,d})$ such that $n_{i} = \sum_{1\leq\alpha\leq r} n_{\alpha,i}$ for every $i=1,\ldots, d$. 
An element of the connected component corresponding to $(\nn_1,\ldots,\nn_r)$ in the right hand side is a tuple of nested quotients 
\begin{align*}
\left( [L_\alpha\twoheadrightarrow T^{(\alpha)}_d\twoheadrightarrow\dots \twoheadrightarrow T^{(\alpha)}_1]\right)_{1\leq \alpha\leq r}, 
\end{align*}
where each $T_i^{(\alpha)}$ is the structure sheaf of a finite subscheme of $X$ of length $n_{\alpha,i}$. By Bifet's theorem on the $\TT$-fixed locus of ordinary Quot schemes \cite{Bifet}, we have that 
\begin{equation}\label{sum_surjections}
\bigoplus_{1\leq \alpha\leq r} \left(L_\alpha\twoheadrightarrow T_i^{(\alpha)}\right) \in \Quot_X(E,n_i)^{\TT}
\end{equation}
for each $i=1,\ldots,d$, and since each of the original tuples of quotients was  nested according to $\nn$, it follows that also the tuples \eqref{sum_surjections} are  nested according to $\nn$, and this proves that \eqref{sum_surjections} defines a point in $\Quot_X(E, \nn)^{\TT}$.

The reverse inclusion follows by an analogous reasoning relying once more on Bifet's result \cite{Bifet}.
\end{proof}

\begin{remark}\label{rmk:rank_one}
For a locally free sheaf $L$ of rank 1, we naturally have the isomorphism
\begin{align*}
     \Quot_X(L, \nn)\cong \Hilb^{\nn}(X),
\end{align*}
where $ \Hilb^{\nn}(X)$ is the nested Hilbert scheme of points, see for example \cite{MR1616606}. Moreover, if $X=C$ is a smooth quasiprojective curve, we have (see \cite[Sec. 0.2]{MR1616606})
\begin{equation}\label{eqn:curve}
     \Hilb^{\nn}(C) \cong \Sym^{n_1}(C)\times \Sym^{n_2-n_1}(C)\times\dots \times \Sym^{n_d-n_{d-1}}(C).
\end{equation}
\end{remark}
Assume now $X=C$ is a smooth quasiprojective curve and let  $x\in \Quot_C(E, \nn)^\TT$ be a $\TT$-fixed point, corresponding to the tuple
\begin{align}\label{eqn: fixed point}
\left( \bigl[L_\alpha\twoheadrightarrow T^{(\alpha)}_d\twoheadrightarrow\dots \twoheadrightarrow T^{(\alpha)}_1\bigr]\right)_\alpha\in \prod_{\alpha=1}^r   \Quot_C(L_\alpha, \nn_\alpha). 
\end{align}
Set $K_i^{(\alpha)}=\ker(L_\alpha\twoheadrightarrow T_i^{(\alpha)})$. The tangent space at $x$ can be written as
\begin{equation}\label{tg_space}
   T_x\Quot_C(E, \nn)=\ker\left(\bigoplus_{1\leq \alpha,\beta\leq r}\bigoplus_{i=1}^d\Hom_C\bigl(K^{(\alpha)}_i, T^{(\beta)}_i\bigr)\xrightarrow{\Delta_x} \bigoplus_{1\leq \alpha,\beta\leq r}\bigoplus_{i=1}^{d-1}\Hom_C\bigl(K^{(\alpha)}_{i+1}, T^{(\beta)}_i\bigr)\right).
\end{equation}
Denote by $w_1,\ldots, w_r$ the coordinates of the algebraic torus $\TT$, which we see as irreducible $\TT$-characters. As a $\TT$-representation, the tangent space  admits the following weight decomposition
\begin{multline*}
   T_x\Quot_C(E, \nn) \\
   =\ker\left(\bigoplus_{1\leq \alpha,\beta\leq r}\bigoplus_{i=1}^d\Hom_C\bigl(K^{(\alpha)}_i\otimes w_\alpha, T^{(\beta)}_i\otimes w_\beta\bigr)\xrightarrow{\Delta_x} \bigoplus_{1\leq \alpha,\beta\leq r}\bigoplus_{i=1}^{d-1}\Hom_C\bigl(K^{(\alpha)}_{i+1}\otimes w_\alpha, T^{(\beta)}_i\otimes w_\beta\bigr)\right).
\end{multline*}
We recall the classical result of Białynicki-Birula (see \cite[Section 4]{BBdecomposition}), by which we obtain a decomposition of $ \Quot_X(E,\nn)$ in the case when $E$ is completely decomposable.

\begin{theorem}[Białynicki-Birula]\label{thm: BB}
Let $X$ be a smooth projective scheme with a $\BG_m$-action and let $\set{X_i}_i$ be the connected components of the $\BG_m$-fixed locus $X^{\BG_m} \subset X$. Then there exists a locally closed stratification $X=\coprod_i X^+_i$, such that each $X^+_i\to X_i$ is an affine fibre bundle. Moreover, for every closed point $x\in X_i$, the tangent space is given by $T_x(X^+_i)=T_x(X)^{\fix}\oplus T_x(X)^{+}$, where $T_x(X)^{\fix}$ (resp. $T_x(X)^{+} $) denotes the $\BG_m$-fixed  (resp. positive) part of $T_x(X)$. In particular, the relative dimension of $X^+_i\to X_i$ is equal to $\dim T_x(X)^{+}$ for $x \in X_i$.
\end{theorem}

The Białynicki-Birula ``strata'' are constructed as follows. If $t$ denotes the coordinate of $\BG_m$, we have
\[
X_i^+=\Set{x\in X \,|\, \lim_{t\to 0}t\cdot x\in X_i}.
\]
In particular, the properness assumption assures that the closure of each $\BG_m$-orbit in $X$ contains the  $\BG_m$-fixed point $\lim_{t\to 0}t\cdot x$. Recently Jelisiejew--Sienkiewicz \cite{MR4021177} generalised Theorem \ref{thm: BB}, proving the the $X^{+}_i$ always exists even when $X$ is not projective (or even not smooth). However, in the smooth case they cover $X$ as long as the closure of every $\BG_m$-orbit contains a fixed point.

\smallbreak
We now determine a  Białynicki-Birula decomposition for $\Quot_C(E, \nn)$, where $C$ is a smooth \emph{quasiprojective} curve. See Mochizuki's paper \cite[Section 2.3.4]{Mochizuki_FILT} for an equivalent construction and tangent space calculation (in the projective case), using a slightly different, but technically equivalent, tangent complex.\footnote{We thank Takuro Mochizuki for kindly sharing with us a note proving that the tangent complex used in \cite{Mochizuki_FILT} is quasi-isomorphic to the one encoded by the map $\Delta_x$.} 

Let $\BG_m\hookrightarrow \TT$ be the generic 1-parameter subtorus given by $w\mapsto (w, w^2,\dots, w^r)$; it is clear that $\Quot_C(E, \nn)^\TT=\Quot_C(E, \nn)^{\BG_m}$. Let
\[
Q_{\underline{\nn}}=\prod_{\alpha=1}^r   \Quot_C(L_\alpha, \nn_\alpha)\subset  \Quot_C(E, \nn)^{\BG_m}
\]
be the connected component of the fixed locus corresponding to the $r$-tuple $\underline{\nn}=(\nn_\alpha)_{1\leq \alpha\leq r}$ decomposing $\nn_1+\dots +\nn_r=\nn$.

\begin{prop}\label{prop: BB for Quot}
Let $C$ be a smooth quasiprojective curve and $E=\bigoplus_{\alpha=1}^rL_\alpha$. Then the nested Quot scheme admits a locally closed stratification 
\[
\Quot_C(E, \nn)=\coprod_{\underline{\nn}}Q^{+}_{\underline{\nn}},
\]
where $\underline{\nn}=(\nn_\alpha)_{1\leq \alpha\leq r}$ are such that  $\nn_1+\dots+\nn_r=\nn$ and $Q^{+}_{\underline{\nn}}\to Q_{\underline{\nn}}$ is an affine fibre bundle of relative dimension $\sum_{1\leq \alpha\leq r}(\alpha-1)n_{\alpha, d}$.
\end{prop}

\begin{proof}
The strata $Q^{+}_{\underline{\nn}}$ are induced by Theorem  \ref{thm: BB} --- we just need to show that the closure of every orbit contains a fixed point.
Choose a compactification $C\hookrightarrow \overline{C}$, an extension $\overline{L}_\alpha$ of each line bundle $L_\alpha$ and consider the induced open immersion 
\[
\Quot_C\left(\bigoplus_{\alpha=1}^rL_\alpha, \nn\right)\hookrightarrow \Quot_{\overline{C}}\left(\bigoplus_{\alpha=1}^r\overline{L}_\alpha, \nn\right).
\]
The closure of every orbit must contain a fixed point in $ \Quot_{\overline{C}}\left(\bigoplus_{\alpha=1}^r\overline{L}_\alpha, \nn\right)$, but the $\BG_m$-action does not move the support of a nested quotient, by which we conclude that such a fixed point had to be  already contained  in $\Quot_C\left(\bigoplus_{\alpha=1}^rL_\alpha, \nn\right)$.

Let $x\in Q_{\underline{\nn}}$ be a fixed point as in \eqref{eqn: fixed point}. The positive part of the tangent space \eqref{tg_space} is
\begin{align*}
   T^+_x\Quot_C(E, \nn)=\ker\left(\bigoplus_{\alpha<\beta}\bigoplus_{i=1}^d\Hom_C\bigl(K^{(\alpha)}_i, T^{(\beta)}_i\bigr)\xrightarrow{\Delta^+_x} \bigoplus_{\alpha<\beta}\bigoplus_{i=1}^{d-1}\Hom_C\bigl(K^{(\alpha)}_{i+1}, T^{(\beta)}_i\bigr)\right),
\end{align*}
where $\Delta^+_x$ is the restriction of the map $\Delta_x$. Thanks to the vanishings \eqref{eqn: vanishings}, $\Delta^+_x$ is surjective, therefore the relative dimension is computed as
\begin{align*}
    \dim_{\bfk}  T^+_x\Quot_C(E, \nn)&=\dim_{\bfk}\left( \bigoplus_{\alpha<\beta}\bigoplus_{i=1}^d\Hom_C\bigl(K^{(\alpha)}_i, T^{(\beta)}_i\bigr)\right)-\dim_{\bfk}\left( \bigoplus_{\alpha<\beta}\bigoplus_{i=1}^{d-1}\Hom_C\bigl(K^{(\alpha)}_{i+1}, T^{(\beta)}_i\bigr)\right)\\
    &=\sum_{\alpha<\beta}\left(\sum_{i=1}^d n_{\beta, i}- \sum_{i=1}^{d-1}n_{\beta, i}\right)\\
    &=\sum_{\beta=1}^r(\beta-1)n_{\beta,d}
\end{align*}
where we used $n_{\beta,i} = \dim_{\bfk} \Hom_C(K_i^{(\alpha)},T_i^{(\beta)})$ since $K_i^{(\alpha)}=\ker(L_\alpha\twoheadrightarrow T_i^{(\alpha)})$ has rank $1$. The proof is complete.
\end{proof}

\section{The motive of the nested Quot scheme on a curve}\label{sec: motives}

\subsection{Grothendieck ring of varieties}\label{sec:K0VAR}
Let $B$ be a scheme locally of finite type over $\bfk$. The \emph{Grothendieck group of $B$-varieties}, denoted $K_0(\Var_{B})$, is defined to be the free abelian group generated by isomorphism classes $[X\to B]$ of finite type $B$-varieties, modulo the scissor relations, namely the identities $[h\colon X \to B] = [h|_Z\colon Z \to B] + [h|_{X\setminus Z} \colon X\setminus Z \to B]$ whenever $Z \into X$ is a closed $B$-subvariety of $X$. The neutral element for the addition operation is the class of the empty variety. The operation
\[
[X\to B]\cdot [X'\to B] = [X\times_{B} X' \to B]
\]
defines a ring structure on $K_0(\Var_{B})$, with identity $\mathbb 1_{B} = [\id \colon B \to B]$. Therefore $K_0(\Var_{B})$ is called the \emph{Grothendieck ring of $B$-varieties}. If $B = \Spec \bfk$, we write $K_0(\Var_{\bfk})$ instead of $K_0(\Var_{\Spec \bfk})$, and we shorten $[X] = [X\to \Spec \bfk]$ for every $\bfk$-variety $X$.

The main rules for calculations in $K_0(\Var_{\bfk})$ are the following:
\begin{enumerate}
    \item If $X \to Y$ is a geometric bijection, i.e.~a bijective morphism, then $[X] = [Y]$.
    \item If $X \to Y$ is Zariski locally trivial with fibre $F$, then $[X] = [Y]\cdot [F]$.
\end{enumerate}

These are, indeed, the only properties that we will use.

The \emph{Lefschetz motive} is the class $\BL = [\BA^1_{\bfk}] \in K_0(\Var_{\bfk})$. It can be used to express, for instance, the class of the projective space, namely $[\BP^n_{\bfk}] = 1+\BL+\cdots+\BL^n \in K_0(\Var_{\bfk})$.

\subsection{Independence of the vector bundle}

The following result generalises \cite[Corollary 2.5]{ricolfi2019motive}, which in turn generalises the main theorem of \cite{BFP19} extending it from proper smooth curves to arbitrary smooth varieties.

\begin{prop}\label{prop: indep vector bundle}
Let $E$ be a locally free sheaf of rank $r$ on a $\bfk$-variety $X$. For every $\nn$, the motivic class of $\Quot_X(E,\nn)$ is independent of $E$, that is
\[
\bigl[\Quot_X(E,\nn)\bigr]=\bigl[\Quot_X(\OO_X^{\oplus r},\nn)\bigr]\in K_0(\Var_{\bfk}).
\]
\end{prop}

\begin{proof}
Let $(U_k)_{1\leq k\leq e}$ be a Zariski open cover trivialising $E$. We can refine it to a locally closed stratification $X = W_1\amalg \cdots\amalg W_e$ such that $W_k \subset U_k$, so that in particular $E|_{W_k} = \OO_{W_k}^{\oplus r}$ for every $k$. Each $W_k$ is taken with the reduced induced scheme structure.

Let $\Quot_{X,W_k}(E,\nn) \subset \Quot_X(E,\nn)$ be the preimage of $\Sym^{n_d}(W_k) \subset \Sym^{n_d}(X)$ along the projection
\[
\pr_d \circ\, \mathsf h_{E,\nn}\colon \Quot_X(E,\nn) \to \prod_{i=1}^d \Sym^{n_i}(X) \to \Sym^{n_d}(X),
\]
where $\mathsf h_{E,\nn}$ is the support map \eqref{supportmap}. We endow $\Quot_{X,W_k}(E,\nn)$ with the reduced scheme structure.
We have a geometric bijection
\[
\coprod_{\nn_1+\cdots+\nn_e = \nn} \prod_{k=1}^e \Quot_{X,W_k}(E,\nn_k) \to \Quot_X(E,\nn),
\]
therefore the motive $[\Quot_X(E,\nn)]$ is computed entirely in terms of the motives $[\Quot_{X,W_k}(E,\nn_k)]$. It is enough to prove that these are independent of $E$.
In the cartesian diagram
\[
\begin{tikzcd}[row sep=large]
\Quot_{U_k,W_k}(E|_{U_k},\nn_k)\MySymb{dr}\arrow[hook]{r}{j}\arrow[hook]{d}
& 
\Quot_{X,W_k}(E,\nn_k)\arrow[hook]{d} \\
\Quot_{U_k}(E|_{U_k},\nn_k)\arrow[hook]{r}{\textrm{open}}
& \Quot_X(E,\nn_k)
\end{tikzcd}
\]
the open immersion $j$ is in fact surjective, hence an isomorphism. But we can repeat this process with $\OO_X^{\oplus r}$ in the place of $E$. It follows that 
\[
    \Quot_{X,W_k}(E,\nn_k) \,\cong\, \Quot_{U_k,W_k}(\OO_{U_k}^{\oplus r},\nn_k)
    \,\cong \,\Quot_{X,W_k}(\OO_X^{\oplus r},\nn_k),
\]
which yields the result.
\end{proof}

\subsection{Proof of the main theorem}\label{sec:MAINTHM_PROOF}
Let $X$ be a smooth quasiprojective variety and $E$ a locally free sheaf of rank $r$. Define
\begin{align*}
    \mathsf{Z}_{X,r,d}(\boldit{q})=\sum_{\nn}\,\bigl[\Quot_X(E, \nn)\bigr]\boldit{q}^{\nn}\in K_0(\Var_{\bfk})\llbracket q_1, \dots, q_d\rrbracket,
\end{align*}
where $\nn=(n_1\leq \dots \leq n_d)$ and we use the multivariable notation $\boldit{q} = (q_1,\ldots,q_d)$ and $\boldit{q}^{\nn}=\prod_{i=1}^dq^{n_i}$. The notation $\mathsf{Z}_{X,r,d}$ reflects the independence on $E$ that we proved in \Cref{prop: indep vector bundle}. If $X=C$ is a smooth quasiprojective curve and $r=d=1$, then $\mathsf{Z}_{C,1,1}(q)$ is simply the Kapranov motivic zeta function
\begin{equation}\label{eqn:kapranov}
\mathsf{Z}_{C, 1,1}(q)=\zeta_C(q)=\sum_{n\geq 0}\,\bigl[\Sym^n(C)\bigr]q^n.
\end{equation}

We can now prove our main theorem, first stated in \Cref{thm:main_thm_INTRO} in the Introduction.

\begin{theorem}\label{mainthm_BODY}
Let $C$ be a smooth quasiprojective curve. The generating series $\mathsf{Z}_{C,r,d}(q)$ is a product of shifted motivic zeta functions: there is an identity
\begin{align*}
     \mathsf{Z}_{C,r,d}(\boldit{q})=\prod_{\alpha=1}^r\prod_{i=1}^d\zeta_C\left(\BL^{\alpha-1}q_iq_{i+1}\cdots q_d\right).
\end{align*}
In particular, $\mathsf{Z}_{C,r,d}(\boldit{q})$ is a rational function in $q_1,\ldots,q_d$.
\end{theorem}
\begin{proof}
By \Cref{prop: indep vector bundle} the motive $[\Quot_C(E,\nn)]$ is independent on the vector bundle $E$, so we may assume $E=\OO_C^{\oplus r}$. In this case, we may compute the motive exploiting the decomposition of $\Quot_C(\OO_C^{\oplus r}, \nn) $ given by \Cref{prop: BB for Quot}. Every stratum is a Zariski locally trivial fibration  over a connected component of the fixed locus, with fibre an affine space whose dimension we computed in \Cref{prop: BB for Quot}.

In what follows, we denote by $\nn_{\alpha}=(n_{\alpha, 1}\leq \dots \leq n_{\alpha, d}) $ a nested tuple of non-negative integers and by $\boldit{l}_{\alpha}=(l_{\alpha,1},\dots,l_{\alpha,d})$ a tuple of non-negative integers. Clearly the two collections of tuples are in bijection, by means of the correspondence 
\begin{equation}\label{correspondence}
(n_{\alpha, 1}\leq \dots \leq n_{\alpha, d}) \longleftrightarrow (n_{\alpha, 1}, n_{\alpha,2}-n_{\alpha,1}, \dots, n_{\alpha, d}-n_{\alpha,d-1}).
\end{equation}
We compute
\begin{align*}
    \sum_{\nn}\,\bigl[\Quot_C(\OO_C^{\oplus r}, \nn) \bigr]\boldit{q}^{\nn}
    &= \sum_{\nn}\boldit{q}^{\nn}\sum_{\nn_1+\dots+\nn_r=\nn}\prod_{\alpha=1}^r\,\bigl[\Quot_C(\OO_C, \nn_\alpha)\bigr]\cdot \BL^{(\alpha-1)n_{\alpha,d}} & \textrm{by \Cref{prop: BB for Quot}}\\
    &=\sum_{\nn_1,\ldots,\nn_r}\prod_{\alpha=1}^r\boldit{q}^{\nn_\alpha}\bigl[\Hilb^{\nn_\alpha}(C)\bigr]\cdot \BL^{(\alpha-1)n_{\alpha,d}}\\
    &= \sum_{\boldit{l}_1,\ldots,\boldit{l}_r}\prod_{\alpha=1}^r\left(\prod_{i=1}^dq_i^{\sum_{j=1}^il_{\alpha,j}}\right)\cdot\bigl[\Hilb^{\nn_\alpha}(C)\bigr]\cdot \BL^{(\alpha-1)\sum_{i=1}^dl_{\alpha,i}} & \textrm{by \eqref{correspondence}}\\
    &=\sum_{\boldit{l}_1,\ldots,\boldit{l}_r}\prod_{\alpha=1}^r\prod_{i=1}^dq_i^{\sum_{j=1}^il_{\alpha,j}}\cdot \bigl[\Sym^{l_{\alpha,i}}(C)\bigr]\cdot \BL^{(\alpha-1)l_{\alpha,i}} & \textrm{by \eqref{eqn:curve}}\\
    &=\sum_{\boldit{l}_1,\ldots,\boldit{l}_r}\prod_{\alpha=1}^r\prod_{i=1}^d\left(q_iq_{i+1}\cdots q_d\right)^{l_{\alpha,i}}\cdot \bigl[\Sym^{l_{\alpha,i}}(C)\bigr]\cdot \BL^{(\alpha-1)l_{\alpha,i}} \\
    &=\prod_{\alpha=1}^r\prod_{i=1}^d\zeta_C\left(\BL^{\alpha-1}q_iq_{i+1}\cdots q_d\right) & \textrm{by \eqref{eqn:kapranov}}.
\end{align*}
The rationality follows by the rationality of the Kapranov zeta function, proved in \cite[Theorem 1.1.9]{Kapranov_rational_zeta}.
\end{proof}

\begin{remark}
We can reformulate our main theorem in terms of the motivic exponential, for which a minimal background is provided in \Cref{appendix}. The case $r=d=1$ yields the classical expression
\[
\zeta_C(q)=\Exp_+([C]q).
\]
The general case becomes
\begin{align*}
\mathsf{Z}_{C,r,d}(\boldit{q}) &= \Exp_+ \left(\bigl[C\bigr] \sum_{\alpha=1}^r \BL^{\alpha-1} \sum_{i=1}^d q_iq_{i+1}\cdots q_d  \right)\\
&=\Exp_+ \left(\bigl[C\times_{\bfk} \BP^{r-1}_{\bfk}\bigr] \sum_{i=1}^d q_iq_{i+1}\cdots q_d  \right).
\end{align*}

Setting $d=1$ we recover the calculations of \cite{BFP19,ricolfi2019motive}.
\end{remark}

\subsection{Hodge--Deligne polynomial}\label{HD_poly}
In this subsection we work over $\bfk = \BC$.
Ring homomorphisms $K_0(\Var_{\BC}) \to R$ are called \emph{motivic measures}. A typical example of a motivic measure is the Hodge--Deligne polynomial
\[
\mathsf E\colon K_0(\Var_{\BC}) \to \BZ[u,v],
\]
defined by sending the class $[Y]$ of a smooth projective variety\footnote{By a beautiful result of Bittner \cite{Bittner01}, the classes of smooth projective varieties generate $K_0(\Var_{\bfk})$ as soon as $\mathrm{char}\, \bfk = 0$. But of course $\mathsf E$ can be defined on arbitrary varieties via mixed Hodge structures.} $Y$ to
\[
\mathsf{E}(Y;u,v) = \sum_{p,q\geq 0}\dim_{\BC} \mathrm{H}^q(Y,\Omega^p_Y)(-u)^p(-v)^q.
\]

\begin{notation}
If $f(u,v) = \sum_{i,j}p_{ij}u^iv^j \in \BZ[u,v]$, we set 
\[
(1-q)^{-f(u,v)} = \prod_{i,j}\,\left(1-u^iv^jq\right)^{-p_{ij}}.
\]
This is actually the formula defining the \emph{power structure} on $\BZ[u,v]$. The motivic measure $\mathsf E$ can be proved to be a morphism of rings with power structure, see \cite{GLMHilb} for full details.
\end{notation}

Let $C$ be a smooth projective curve of genus $g$. We have
\begin{equation} \label{E_curve}
\begin{split}
\mathsf E(\zeta_C(q)) &= \sum_{n\geq 0} \mathsf E(\Sym^n(C);u,v)q^n = (1-q)^{-\mathsf E(C;u,v)} \\
& =  (1-q)^{-(1-gu-gv+uv)} \\
& =\frac{(1-uq)^g(1-vq)^g}{(1-q)(1-uvq)}.    
\end{split}
\end{equation}
For $E$ a locally free sheaf of rank $r$ over $C$, define
\[
\mathsf E_{C,r,d}(\boldit{q}) = \sum_{\nn} \mathsf E(\Quot_C(E,\nn);u,v)\boldit{q}^{\nn}.
\]
As a direct consequence of \Cref{mainthm_BODY}, we obtain the following corollary.

\begin{corollary}
There is an identity
\[
\mathsf E_{C,r,d}(\boldit{q}) 
= \prod_{\alpha=1}^r\prod_{i=1}^d \frac{\left(1-u^\alpha v^{\alpha-1}q_iq_{i+1}\cdots q_d\right)^g\left(1-u^{\alpha-1}v^\alpha q_iq_{i+1}\cdots q_d\right)^g}{\left(1-u^{\alpha-1}v^{\alpha-1}q_iq_{i+1}\cdots q_d\right)\left(1-u^{\alpha}v^{\alpha}q_iq_{i+1}\cdots q_d\right)}.
\]
\end{corollary}

\begin{proof}
This follows by combining \Cref{mainthm_BODY} and \Cref{E_curve} with one another, after observing that $\mathsf E$ is multiplicative (being a ring homomorphism) and sends $\BL$ to $uv$.
\end{proof}

The generating function of the signed Poincar\'e polynomials is obtained from $\mathsf E_{C,r,d}(\boldit{q})$ by setting $u=v$. The result confirms a result of L.~Chen \cite{Chen-hyperquot} obtained in the case $C=\BP^1$. The generating series of topological Euler characteristics is obtained from $\mathsf E_{C,r,d}(\boldit{q})$ by setting $u=v=1$, also in the quasiprojective case. So we obtain
\[
\sum_{\nn} e_{\mathrm{top}}(\Quot_C(E,\nn)) \boldit{q}^{\nn} 
= \prod_{i=1}^d\left(1-q_iq_{i+1}\cdots q_d\right)^{-r\cdot e_{\mathrm{top}}(C)}.
\]

\appendix
\section{Motivic exponentials}\label{appendix}

If $(\Lambda,\mu,\epsilon)$ is a commutative monoid in the category of $\bfk$-schemes, where $\mu \colon \Lambda \times \Lambda \to \Lambda$ is the multiplication map and $\epsilon\colon \Spec \bfk \to \Lambda$ is the identity element, then by \cite[Example 3.5\,(4)]{BenSven3}, one has a $\lambda$-ring structure $\sigma_\mu$ on the Grothendieck ring
\[
K_0(\Var_{\Lambda}),
\]
determined by the operations
\[
\sigma_\mu^n\,\bigl[Y \xrightarrow{f} \Lambda\bigr] \,=\, \bigl[\Sym^n Y \xrightarrow{\Sym^nf} \Sym^n\Lambda \xrightarrow{\mu} \Lambda\bigr].
\]
Assume $\Lambda_+ \subset \Lambda$ is a sub-monoid such that $\coprod_{n>0}\Lambda_+^{\times n} \to \Lambda$ is of finite type. Then we can define the \emph{motivic exponential}
\[
\Exp_\mu\colon K_0(\Var_{\Lambda_+}) \to K_0(\Var_{\Lambda})^\times
\]
by setting
\[
\Exp_\mu(A) = \sum_{n\geq0}\sigma^n_\mu(A)
\]
for an effective class $A$, and extending via
\[
\Exp_\mu(A-B) = \Exp_\mu(A)\cdot \Exp_\mu(B)^{-1}
\]
whenever $A$ and $B$ are effective. The map $\Exp_\mu$ is injective.
See \cite[Section 1]{DavisonR} for more details.

We use this construction in the case $(\Lambda,\mu,\epsilon) = (\BN^d,+,0)$, and setting $\Lambda_+ = \BN^d\setminus 0$. Of course here we are seeing $\BN^d$ as the $\bfk$-scheme $\coprod_{\nn \in \BN^d}\Spec \bfk$. There is an isomorphism
\[
\begin{tikzcd}
K_0(\Var_{\bfk})\llbracket q_1,\ldots,q_d \rrbracket  \arrow{r}{\sim} & K_0(\Var_{\BN^d})
\end{tikzcd}
\]
defined by sending
\[
\sum_{\nn \in \BN^d} Y_{\nn} \cdot q_1^{n_1}\cdots q_d^{n_d} \,\,\mapsto \,\, \left[\coprod_{\nn \in \BN^d} Y_{\nn} \to \Spec \bfk(\nn) \right] 
\]
for varieties $Y_{\nn}$, and extending by linearity. Under this identification, if we let $\mathfrak m$ be the ideal generated by $(q_1,\ldots,q_d)$ in $K_0(\Var_{\bfk}) \llbracket q_1,\ldots,q_d \rrbracket$, we can see $\Exp_+$ as a group isomorphism
\[
\begin{tikzcd}
\Exp_+\colon \mathfrak m\cdot K_0(\Var_{\bfk})\llbracket q_1,\ldots,q_d \rrbracket \arrow{r}{\sim} &
1+ \mathfrak m \cdot K_0(\Var_{\bfk}) \llbracket q_1,\ldots,q_d \rrbracket \subset (K_0(\Var_{\bfk}) \llbracket q_1,\ldots,q_d\rrbracket)^\times
\end{tikzcd}
\]
between an additive group (on the left) and a multiplicative group (on the right). In particular, one has the identity
\[
\Exp_+\left(\sum_{\ell=1}^s f_\ell(q_1,\ldots,q_d)\right) = \prod_{\ell=1}^s \Exp_+(f_\ell(q_1,\ldots,q_d))
\]
for $f_\ell(q_1,\ldots,q_d) \in \mathfrak m\cdot K_0(\Var_{\bfk})\llbracket q_1,\ldots,q_d \rrbracket$.

\bibliographystyle{amsplain-nodash} 
\bibliography{bib}

\ifx\undefined\bysame
\newcommand{\bysame}{\leavevmode\hbox to3em{\hrulefill}\,}
\fi
\begin{thebibliography}{10}

\bibitem{BFP19}
Massimo Bagnarol, Barbara Fantechi, and Fabio Perroni, {\em {On the motive of
  zero-dimensional Quot schemes on a curve}}, New York J. Math. {\bf 26}
  (2019), no.~2020, 138--148.

\bibitem{BR18}
Sjoerd Beentjes and Andrea~T. Ricolfi, {\em {Virtual counts on Quot schemes and
  the higher rank local DT/PT correspondence}}, {Math. Res. Lett.} {\bf 28}
  (2021), no.~4, 967--1032.

\bibitem{BBdecomposition}
Andrzej Bia{\l}ynicki-Birula, {\em Some theorems on actions of algebraic
  groups}, Ann. of Math. (2) {\bf 98} (1973), 480--497.

\bibitem{Bifet}
Emili Bifet, {\em {Sur les points fixes du sch\'{e}ma
  {${\textrm{Quot}}_{{\mathscr O}^r_X/X/k}$} sous l'action du tore {${\mathbf
  G}^r_{m,k}$}}}, C. R. Acad. Sci. Paris S\'{e}r. I Math. {\bf 309} (1989),
  no.~9, 609--612.

\bibitem{Bittner01}
Franziska {Bittner}, {\em {The universal Euler characteristic for varieties of
  characteristic zero}}, {Compos. Math.} {\bf 140} (2004), no.~4, 1011--1032
  (English).

\bibitem{cazzaniga2020higher}
Alberto Cazzaniga, Dimbinaina Ralaivaosaona, and Andrea~T. Ricolfi, {\em
  {Higher rank motivic Donaldson--Thomas invariants of $\mathbb{A}^3$ via
  wall-crossing, and asymptotics}}, Math. Proc. Cambridge Philos. Soc.,
  \href{https://www.cambridge.org/core/journals/mathematical-proceedings-of-the-cambridge-philosophical-society/article/higher-rank-motivic-donaldsonthomas-invariants-of-mathbb-a3-via-wallcrossing-and-asymptotics/4181C5A9585E17889AD9CB678784D55B}{Online
  first}, 2022.

\bibitem{Cazzaniga:2020aa}
Alberto Cazzaniga and Andrea~T. Ricolfi, {\em {Framed motivic Donaldson--Thomas
  invariants of small crepant resolutions}}, Math. Nachr. {\bf 295} (2022),
  no.~6, 1096--1112.

\bibitem{MR2716513}
Jan Cheah, {\em The cohomology of smooth nested {H}ilbert schemes of points},
  ProQuest LLC, Ann Arbor, MI, 1994, PhD Thesis. The University of Chicago.

\bibitem{MR1616606}
Jan Cheah, {\em Cellular decompositions for nested {H}ilbert schemes of
  points}, Pacific J. Math. {\bf 183} (1998), no.~1, 39--90.

\bibitem{MR1612677}
Jan Cheah, {\em The virtual {H}odge polynomials of nested {H}ilbert schemes and
  related varieties}, Math. Z. {\bf 227} (1998), no.~3, 479--504.

\bibitem{Chen-hyperquot}
Linda Chen, {\em Poincar{\'e} polynomials of hyperquot schemes}, Math. Ann.
  {\bf 321} (2001), no.~2, 235--251.

\bibitem{BenSven3}
Ben Davison and Sven Meinhardt, {\em Motivic {D}onaldson--{T}homas invariants
  for the one-loop quiver with potential}, Geom. Topol. {\bf 19} (2015), no.~5,
  2535--2555.

\bibitem{DavisonR}
Ben Davison and Andrea~T. Ricolfi, {\em {The local motivic DT/PT
  correspondence}}, J. Lond. Math. Soc. {\bf 104} (2021), no.~3, 1384--1432.

\bibitem{dZNPZ_playing_index_M_theory}
Michele Del~Zotto, Nikita Nekrasov, Nicolò Piazzalunga, and Maxim Zabzine,
  {\em {Playing with the index of M-theory}},
  {\href{https://arxiv.org/abs/2103.10271}{ArXiv:2103.10271}} 2021.

\bibitem{FMR_K-DT}
Nadir Fasola, Sergej Monavari, and Andrea~T. Ricolfi, {\em {Higher rank
  K-theoretic Donaldson--Thomas theory of points}}, Forum Math. Sigma {\bf 9}
  (2021), no.~E15, 1--51.

\bibitem{Fogarty_Hilb}
John Fogarty, {\em Algebraic families on an algebraic surface}, Amer. J. Math.
  {\bf 90} (1968), 511--521.

\bibitem{Heinloth_altri}
Oscar Garc\'{\i}a-Prada, Jochen Heinloth, and Alexander Schmitt, {\em On the
  motives of moduli of chains and {H}iggs bundles}, J. Eur. Math. Soc. (JEMS)
  {\bf 16} (2014), no.~12, 2617--2668.

\bibitem{GLMHilb}
Sabir~M. {Gusein-Zade}, Ignacio {Luengo}, and Alejandro {Melle-Hern\'andez},
  {\em {Power structure over the Grothendieck ring of varieties and generating
  series of Hilbert schemes of points}}, {Mich. Math. J.} {\bf 54} (2006),
  no.~2, 353--359.

\bibitem{MR0213368}
Robin Hartshorne, {\em Connectedness of the {H}ilbert scheme}, Inst. Hautes
  \'{E}tudes Sci. Publ. Math. (1966), no.~29, 5--48.

\bibitem{hoskins2019formula}
Victoria Hoskins and Simon Pepin~Lehalleur, {\em A formula for the {V}oevodsky
  motive of the moduli stack of vector bundles on a curve}, Geom. Topol. {\bf
  25} (2021), no.~7, 3555--3589.

\bibitem{modulisheaves}
Daniel Huybrechts and Manfred Lehn, {\em The geometry of moduli spaces of
  sheaves}, Second Edition, Cambridge University Press, 2010, pp.~xviii+325.

\bibitem{MR4021177}
Joachim Jelisiejew and {\L}ukasz Sienkiewicz, {\em Bia\l ynicki-{B}irula
  decomposition for reductive groups}, J. Math. Pures Appl. (9) {\bf 131}
  (2019), 290--325.

\bibitem{Kapranov_rational_zeta}
Mikhail Kapranov, {\em {The elliptic curve in the S-duality theory and
  Eisenstein series for Kac-Moody groups}},
  {\href{https://arxiv.org/abs/math/0001005}{ArXiv:0001005}}, 2000.

\bibitem{Mochizuki_FILT}
Takuro Mochizuki, {\em {The structure of the cohomology ring of the filt
  schemes}}, {\href{https://arxiv.org/abs/math/0301184}{ArXiv:0301184}}, 2003.

\bibitem{Lissite-quot}
Sergej Monavari and Andrea~T. Ricolfi, {\em {Sur la lissit\'e du sch\'ema Quot
  ponctuel embo\^{i}t\'e}}, Canad. Math. Bull.,
  \href{https://www.cambridge.org/core/journals/canadian-mathematical-bulletin/article/sur-la-lissite-du-schema-quot-ponctuel-emboite/7DEE062EC8DCE202FFC47CEB35CDFD6B}{Online
  first} (2022).

\bibitem{Magnificent_colors}
Nikita Nekrasov and Nicol{\`o} Piazzalunga, {\em Magnificent four with colors},
  Comm. Math. Phys. {\bf 372} (2019), no.~2, 573--597.

\bibitem{ricolfi2019motive}
Andrea~T. Ricolfi, {\em {On the motive of the Quot scheme of finite quotients
  of a locally free sheaf}}, J. Math. Pures Appl. {\bf 144} (2020), 50--68.

\bibitem{Virtual_Quot}
Andrea~T. Ricolfi, {\em {Virtual classes and virtual motives of Quot schemes on
  threefolds}}, Adv. Math. {\bf 369} (2020), 107182.

\bibitem{Sernesi}
Edoardo Sernesi, {\em Deformations of algebraic schemes}, Grundlehren der
  Mathematischen Wissenschaften [Fundamental Principles of Mathematical
  Sciences], vol. 334, Springer-Verlag, Berlin, 2006.

\end{thebibliography}
\end{document}